\documentclass[12pt,reqno]{amsart}

\usepackage[a4paper,top=3.1cm,bottom=3.1cm,left=2.8cm,right=2.8cm,marginparwidth=54pt]{geometry}

\usepackage{latexsym}
\usepackage{amsthm, amsmath, amssymb,latexsym,amsfonts,verbatim}
\usepackage{mathrsfs}

\usepackage{xcolor}
\usepackage{graphicx}
\usepackage{caption}

\usepackage{paralist}
\usepackage[mathscr]{eucal}
\usepackage[pdftex]{hyperref}
\usepackage[capitalise]{cleveref}
\hypersetup{citecolor=blue,linktocpage}
\usepackage[colorinlistoftodos]{todonotes}
\usepackage[all]{xy}

\newcommand{\rk}{\mbox{rank}}

\xyoption{all} \numberwithin{equation}{section}

\newtheorem{theorem}{Theorem}[section]
\newtheorem{lemma}[theorem]{Lemma}
\newtheorem{corollary}[theorem]{Corollary}
\newtheorem{proposition}[theorem]{Proposition}
\theoremstyle{definition}
\newtheorem{remark}[theorem]{Remark}
\newtheorem{definition}[theorem]{Definition}
\newtheorem{example}[theorem]{Example}

\newtheorem{question}[theorem]{Question}

\newif\ifpdf\ifx\pdfoutput\undefined\pdffalse\else\pdfoutput=1\pdftrue\fi

\newcommand{\nc}{\newcommand} 
\nc{\cH}{{\mathcal H}}
\nc{\cA}{{\mathcal A}}
\nc{\cG}{{\mathcal G}}
\nc{\cC}{{\mathcal C}}
\nc{\cO}{{\mathcal O}}
\nc{\cI}{{\mathcal I}}
\nc{\cR}{{\mathcal R}}
\nc{\cB}{{\mathcal B}}
\nc{\cY}{{\mathcal Y}}
\nc{\cK}{{\mathcal K}} 
\nc{\cX}{{\mathcal X}}
\nc{\cS}{{\mathcal S}}
\nc{\cE}{{\mathcal E}}
\nc{\cF}{{\mathcal F}}
\nc{\cZ}{{\mathcal Z}}
\nc{\cQ}{{\mathcal Q}}
\nc{\cN}{{\mathcal N}}
\nc{\cP}{{\mathcal P}}
\nc{\cL}{{\mathcal L}}
\nc{\cM}{{\mathcal M}}
\nc{\cT}{{\mathcal T}}
\nc{\cW}{{\mathcal W}}
\nc{\cU}{{\mathcal U}}
\nc{\cJ}{{\mathcal J}}
\nc{\cV}{{\mathcal V}}
\nc{\bH}{{\mathbb H}}
\nc{\bA}{{\mathbb A}}
\nc{\bG}{{\mathbb G}}
\nc{\bC}{{\mathbb C}}
\nc{\bO}{{\mathbb O}}
\nc{\bI}{{\mathbb I}}
\nc{\bB}{{\mathbb B}}
\nc{\bY}{{\mathbb Y}}
\nc{\bK}{{\mathbb K}} 
\nc{\bX}{{\mathbb X}}
\nc{\bS}{{\mathbb S}}
\nc{\bE}{{\mathbb E}}
\nc{\bF}{{\mathbb F}}
\nc{\bZ}{{\mathbb Z}}
\nc{\bQ}{{\mathbb Q}}
\nc{\bN}{{\mathbb N}}
\nc{\bP}{{\mathbb P}}
\nc{\bL}{{\mathbb L}}
\nc{\bM}{{\mathbb M}}
\nc{\bT}{{\mathbb T}}
\nc{\bW}{{\mathbb W}}
\nc{\bU}{{\mathbb U}}
\nc{\bD}{{\mathbb D}}
\nc{\bJ}{{\mathbb J}}
\nc{\bV}{{\mathbb V}}
\nc{\bbZ}{{\mathbb Z}}
\nc{\bR}{{\mathbb R}}
\nc{\fr}{{\rightarrow}}
\nc{\co}{{\nabla}}

\nc{\cu}{{\barline{\nabla}}}
\newcommand{\pr}{\mathbb P}
\newcommand{\pic}{\mbox{\upshape{Pic}}}
\nc{\divi}{\mbox{\upshape{div}}}

\begin{document}
\begin{abstract}
We study pencils of plane cubics with only one base point and general member smooth, giving a complete classification. Under the additional hypothesis that all members are irreducible, we prove that there exists a unique non-isotrivial pencil with these properties up to projective transformation. 
We compare our construction with the classical approaches given by Gattazzo, Beauville and Miranda-Persson.
\end{abstract}


\title {Pencils of plane cubics with one base point} 
\author{R.~Moschetti, G.P.~Pirola, L.~Stoppino}

\maketitle


\section{Introduction}\label{intro}

The main question from which this work started is the following:
\begin{question}\label{quest: main} Let $d\geq 3$ be any integer. Do there exist pencils of curves of degree $d$ in $\pr^2$ that satisfy the following three conditions?
\begin{enumerate}
\item[(C1)] there is only one base point;
\item[(C2)] the general member is smooth;
\item[(C3)] all the members are irreducible.
\end{enumerate}
\end{question}

If we drop the irreducibility request in \Cref{quest: main}, it is easy to construct such pencils, see \Cref{rem: esempio ovvio}. However, condition (C3) has a more subtle answer. Once we deal with the existence part, we can study deeper properties of such pencils.

\begin{question}\label{quest: canonicity}
Assume \Cref{quest: main} has a positive answer.
\begin{itemize}
\item[(i)] 
What is the codimension in the Severi variety of the locus of plane curves admitting such pencils? 
\item[(ii)] 
Is there a canonical way to obtain these pencils from a suitable curve?
\item[(iii)] Are these pencils non-isotrivial? 
\item[(iv)] Can we classify the singularities of the members of such pencils?
\end{itemize}
\end{question}

The main result of this paper gives a positive answer to \Cref{quest: main}, and addresses all the points of \Cref{quest: canonicity} for the case of plane curves of degree $3$. We work over the complex field $\mathbb C$.

\begin{theorem}\label{teo: mix}
Every smooth cubic $C\subset \pr^2$ belongs to a pencil   $|\cV|$ satisfying all the properties of \Cref{quest: main}. 
Moreover, $|\cV|$ is unique up to projective transformations. The pencil has four singular members, which are all nodal, therefore it is non-isotrivial.
\end{theorem}
Thanks to this result, in this paper we will fix a realisation of such pencil and always denote it by $|\cV|$. In \Cref{ssec: costruzione} and \ref{ssec: unicita'} we give the main construction for pencils satisfying all the properties of \Cref{quest: main}. In \Cref{sec: singular} we study their singular elements. Finally, uniqueness up to projective transformations of $\bP^2$ is proved in \Cref{ssec: uniqueness}.

The key point of the construction of $|\cV|$ is the fact that the base point necessarily is a $9$-torsion point which is not $3$-torsion of a smooth member (once the origin has been fixed in a flex point). 
The theory of points of order nine on a cubic has a long and glorious history, dating back more than 150 years with the works of Cayley \cite{Cayley}, Fano \cite{Fano}, Hart \cite{Hart}, Halphen \cite{Halphen} and  Salmon \cite{Salmon}. More recently, the same kind of topics were considered by Remo Gattazzo in \cite{Gat}, together with a modern viewpoint on the beautiful construction of the so-called tangential triangles of a cubic.\footnote{\emph{``It may be proposed to determine the points at which it shall be possible to draw cubics osculating the given curve in nine points; or in other words, such that the third tangential D may coincide with the original point A.''} Salmon, \cite{Salmon}. See also \cite[24 pag.30]{SK}.} To the present days, the problem is considered for instance in \cite{BB} and in \cite{MZ}. In the recent preprint \cite{Kollar} the same pencil is constructed in Example 46. We give in \Cref{ssec: gattazzo} an overview of Gattazzo's results in relation to ours. Using both these techniques and our approach we can produce explicit examples of equations for $|\cV|$: see \Cref{ex: Ric} and \ref{ex: Gatt}.

The points of $3$-torsion on a plane cubic are also worth studying: they are related with the pencils satisfying (C1), (C2) but not (C3), see \Cref{prop:classification}. We describe completely these pencils and their singular fibres in \Cref{prop: flex}. 
Interestingly enough, some of them are isotrivial, see \Cref{rem: isotriviality}. 
\medskip

Pencils of cubics and elliptic fibrations have been the object of a thorough study over the years. 
We relate our construction with two other perspectives, as follows. 
First of all, the fibration associated with $|\cV|$, constructed by blowing up $9$ infinitely close points, is a non-isotrivial semistable genus $1$ fibration over $\bP^1$ with the minimal number of singular fibres, and so it belongs to Beauville's classification \cite{beau1}. In \Cref{ssec: beauville}  we prove that this fibration is the modular family associated to the subgroup $\Gamma_0(9) \cap \Gamma^0_0(3)$.
We also describe an explicit Cremona transformation of degree $4$ sending $|\cV|$ into the pencil given by Beauville (\Cref{prop: cremona 4}).

In Section \ref{sec: Miranda-Persson} we see that all fibrations associated to pencils satisfying (C1) and (C2) are rational extremal elliptic fibrations (see \Cref{def: reef}), and we show where they appear in Miranda-Persson's classification \cite[Theorem 4.1]{MP}. \medskip

The study carried out in this paper, and more in generally answers to \Cref{quest: main} for any degree, are related with many other research fields. For instance, in \cite[Example 3.4]{BaP} a similar reasoning is used to construct semi-canonical points on curves. Moreover, pencils satisfying (C1) and (C2) are degree $d$ integrable foliations on the projective plane minus one point. Finally, let us note that in the recent paper \cite{weichen}, Wei Chen uses the pencil studied here  to give an example related to the hyperbolicity of the Hirzebruch surface $\mathbb{F}_1$ minus a curve $B$ such that $(\mathbb{F}_1,B)$ is of log general type. 


In a forthcoming paper we study \Cref{quest: main} for the case $d=4$, proving in particular that the general quartic plane curve does not belong to a pencil satisfying all the hypotheses, but that examples exist.

\subsection*{Acknowledgements}  We are grateful to Alice Garbagnati, for explaining us the classification of Miranda and Persson. We owe her \Cref{sec: Miranda-Persson}. We also thank Claudia Alc\`antara, Ciro Ciliberto and Letterio Gatto for their interest in our work and their kind encouragement. We are pleased to thank J\'anos Koll\'ar for suggesting  \cref{rem:kollardue} related to his recent preprint \cite{Kollar}, and giving us some interesting references \cref{rem: kollar}. We are grateful to Wei Chen, for finding some minor mistakes in the previous version of the paper and kindly point them out to us. Finally, thank an anonymous referee for spotting some inaccuracies in a previous version of the paper. 

R.M. is partially supported by the PRIN project 2022L34E7W ``Moduli Spaces and Birational Geometry''. G.P.P. and L.S. are partially supported by the PRIN project 20228JRCYB ``Moduli spaces and special varieties''. All the authors are members of the GNSAGA - INdAM.


\section{Construction of the pencil and an explicit example}\label{ssec: costruzione}
As a first example, we give a straightforward construction of pencils of plane curves of any degree $d\geq 3$ satisfying only (C1) and (C2).

\begin{example}\label{rem: esempio ovvio}
 Let $C\colon F=0$ be a smooth degree $d$ curve possessing a hyperflex $p\in C$, i.e. a flex of maximal order $d$, let $l\colon L=0$ be the tangent line to $C$ at $p$. Then $tF+uL^d=0$, $(t:u)\in \pr^1$ represents a pencil with only $p$ as base point, general smooth member and at least one non-reduced member: $dl$.
 Note that the locus of plane curve possessing a hyperflex is of codimension $d-3$ in the Severi variety of plane curves of degree $d$. 
 Interestingly, in some cases these pencils are isotrivial, see \Cref{rem: isotriviality}.
\end{example}

We now describe pencils $|V|$ of cubics satisfying all the properties of \Cref{quest: main}. We start from any smooth cubic $C\subset \pr^2$ with equation $ F=0$. 
Let us fix one of the $9$ flexes of $C$ as the origin $O$ of the group law on $C$. This identifies the $3$-torsion points with the flexes.
\begin{proposition}\label{prop: construction}
With the above notations, for any point $p$ of $9$-torsion which is not a flex of $C$, there exists an unique pencil containing $C$ satisfying all assumptions of \Cref{quest: main} having $p$ as base point.
\end{proposition}
\begin{proof}
From the assumption on $p$ we have that $\cO_C(9p)\cong \cO_C(9O)$ but $\cO_C(3p)\not\cong \cO_C(3O)$. 
Let us consider the sequence defining $C$:
\[
0\longrightarrow \cO_{\pr^2}(-C)\longrightarrow \cO_{\pr^2}\longrightarrow \cO_C\longrightarrow 0
\]
and tensor it with $\cO_{\pr^2}(3)$:
\[
0\longrightarrow \cO_{\pr^2}\longrightarrow \cO_{\pr^2}(3)\longrightarrow \cO_C\otimes\cO_{\pr^2}(3)\longrightarrow 0.
\]
Observe that $ \cO_C\otimes\cO_{\pr^2}(3)\cong \cO_C(9O)\cong \cO_C(9p)$: indeed, let $H\subset \pr^2$ be the tangent line to $C$ at $O$;. We have 
\[ \cO_C\otimes\cO_{\pr^2}(3)\cong\cO_C(C\cdot 3H)\cong \cO_C(9O).\]

Let us consider the associated cohomology sequence 
\[
0\longrightarrow H^0(\pr^2, \cO_{\pr^2}) \longrightarrow H^0(\pr^2, \cO_{\pr^2}(3))\stackrel{\alpha}\longrightarrow H^0(C,\cO_C(9p)) \longrightarrow  H^1(\pr^2, \cO_{\pr^2})=0.
\]
Let $W\subset H^0(C,\cO_C(9p))$ be the one dimensional subspace of sections with support exactly $9p$. Clearly $\dim \alpha^{-1}(W)=2$, so the pencil that we seek is precisely $V:=\alpha^{-1}(W)$. The general element of $|V|$ is smooth, because $C$ itself is smooth.
So, choosing $\sigma \in H^0(C,\cO_C(9p))$ generating $W$, the general $G\in H^0(\pr^2, \cO_{\pr^2}(3))$ whose image is $\sigma$ is smooth. 

The equation of any member of the pencil is thus generated by $F$ and $G$. Let us prove that this pencil meets the conditions of \Cref{quest: main}. Clearly the only base point is $p$. 

Let us suppose now by contradiction  that there exists a reducible member $D=H+Q\in |V|$, where $H$ is a line, and $Q$ is a (possibly degenerate) conic. 
The local intersection of $D $ and $C$ at $p$ is:
$9=(D\cdot C)_p= (H\cdot C)_p+ (Q\cdot C)_p$.
Now, $(H\cdot C)_p\leq 3$ and $(Q\cdot C)_p\leq 6$, so both are equalities. 
This implies that $H$ is the tangent line to $C$ at $p$, and that $p$ is a flex for $C$, contrary to the assumption.
\end{proof}
\begin{example}\label{ex: Ric}
In order to compute an explicit example of such pencils, start with the nodal cubic $C$ given by $x^3 - y^3 - xyz=0$, which has flexes at the points $(0:0:1)$, $(1:1:0)$, $(1:\rho:0)$, $(1:\rho^2:0)$, where $\rho$ is a cubic root of $1$. In the affine plane given by $y=0$, the cubic becomes $x^3-xz-1=0$, and the projection $\pi_x$ to the $x$ axis gives a morphism to $\bC^*$ compatible with the group law on $C$. We then fix a $9$-th root of $1$ in $\bC^*$, denote it by $\eta$, and such that $\eta^3 \neq 1$. We use $\pi_x$ to get the corresponding point of order $9$ on the curve $C$: $(\eta, (\eta^3-1)/\eta)$. Finally, we can find the family of cubics by using the fact that they all intersects in $(\eta, (\eta^3-1)/\eta)$, starting from a general equation and imposing linear conditions on the coefficients. We get the following pencil:
\begin{align*}
    \big[(-84\eta^3+3)x^3+(-9\eta^8+36\eta^2)x^2y-9\eta xy^2+y^3+(-36\eta^7+126\eta^4-9\eta)x^2+\\
    9\eta^8y^2+(9\eta^8-126\eta^5+36\eta^2)x+(36\eta^7-9\eta)y+84\eta^6-3\big]-t\big[x^3 - xy - 1\big].
\end{align*}
Note that the cubic $C$ is one of the four singular members of the pencil.
\end{example}

\begin{remark}\label{rem: sex}
Recall that for any smooth curve $C$ with degree $d\geq 3$ and for any $p\in C$ there exists a unique conic $Q$ such that 
\[
(C\cdot Q)_p\geq 5=h^0(\pr^2, \cO_{\pr^2}(2)),
\]
called the {\em osculating conic to $C$ in $p$}. In the same way, for any curve $C$ of degree $d\geq 4$ and for any $p\in C$, there exists a unique cubic $D$, which we call  {\em osculating cubic in $p$}, such that 
\[
(C\cdot D)_p\geq 9=h^0(\pr^2, \cO_{\pr^2}(3)).
\]

If we try to extend the definition to lower degrees, we have that the osculating conic of a conic in  any point is the conic itself, hence unique. On the other hand, the osculating cubic of a cubic in a given point $p$ is not necessarily unique: our construction proves indeed that for a smooth cubic, if the point $p$ is of order nine, then there is a 1-parameter family of irreducible cubics intersecting $C$ in $9p$.
\end{remark}

\section{Uniqueness of the construction}\label{ssec: unicita'}
We have seen that any degree $3$ curve can be an element of a pencil satisfying the assumptions of \Cref{quest: main}, thus answering to point  $(i)$ of \Cref{quest: canonicity}: the codimension is in this case $0$.
We now see that any pencil of degree $3$ satisfying the conditions of \Cref{quest: main} can be realised with the construction of \Cref{prop: construction}. Let us first recall the following elementary result on the local intersection of plane curves, which we will repeatedly use in the paper. For a proof see for instance \cite[Lemma 1.3.8]{namba}.
\begin{lemma}\label{lem: namba}
Let $C,D$ and $E\subset \pr^2$ be curves such that $D$ is non-singular in $p$. Then
\[(C\cdot E)_p\geq \min \{(D\cdot E)_p, (D\cdot C)_p\}.\]
\end{lemma}

\begin{proposition}\label{prop: unicita'}
Let $V\subset H^0(\pr^2,\cO_{\pr^2}(3))$ be a pencil of cubics in $\pr^2$ satisfying the conditions of \Cref{quest: main}. Let $p\in \pr^2$ be the base point of $|V|$. Then, for any smooth member  $C\in |V|$, if we fix the origin of the group law on $C$ in a flex, we have that $p$ is a $9$-torsion point of $C$ which is not $3$-torsion. 
\end{proposition}
\begin{proof}
We fix the origin $O$ in a flex of $C$ different than $p$ (we will prove that $p$ is not a flex indeed). Given any other smooth member $D\in |V|$, we have $(C\cdot D)=9p$, so $9p\sim 9O$ in $\pic^9(C)\cong C$, so $p$ is $9$-torsion. 
Suppose now by contradiction that $p$ is a flex for $C$. Let $L$ be the tangent line to $C$ at $p$. Note that $p$ has to be a flex for $D$. Indeed, we have by \Cref{lem: namba} that 
$$(D\cdot L)_p\geq \min \{(C\cdot L)_p, (C\cdot D)_p\}=3.$$
Clearly, no element of  $|V|$ can have $L$ as a component by (C3). Observe now that for any $q\in \pr^2\setminus \{p\}$ there exists a unique $D_q\in |V|$ such that $q\in D_q$. So, let us choose a general point $q\not= p$ in $L$. 
Then, we have $(D_q\cdot L)=q+3p$, which is a contradiction, because $D_q$ has degree $3$.
\end{proof}

\begin{proposition} \label{prop:classification}
Let $V\subset H^0(\pr^2,\cO_{\pr^2}(3))$ be a pencil of cubics in $\pr^2$ satisfying (C1) and (C2). Then either $V$ is as in \Cref{prop: unicita'}, or the base point $p$ of $V$ is a flex for any smooth element $C\in |V|$ and the pencil is generated by $C$ and $3L$ where $L $ is the tangent line to $C$ at $p$ (in other words, $|V|$ is as in \Cref{rem: esempio ovvio}).
\end{proposition}
\begin{proof}
We have seen in Proposition \ref{prop: unicita'} that any pencil of cubics satisfying (C1) and (C2)  has the base point $p$ which is a $9$-torsion point for any smooth element $C\in |V|$. 
If it is not $3$-torsion, (C3) holds. If, on the other hand, $p$ is $3$-torsion for $C$, then necessarily $L$ is in the support of some element of $|V|$; 
indeed, if this was not the case, for any point $q\in L\setminus \{p\}$, there should be an element of the pencil $D_q$ passing through $q$. This leads to the same contradiction as in the proof of \Cref{prop: unicita'}. 
Therefore  there is an element of $|V|$ of the form $L+Q$, where $Q$ is a conic such that $(Q\cdot C)=6p$. 
We see now that  necessarily $Q=2L$. 
Indeed, if $Q$ were irreducible (hence smooth) we would have $(Q\cdot L)_p\geq \min \{(Q\cdot C)_p, (L\cdot C)_p\}=3$ by \Cref{lem: namba}, and this is impossible.
\end{proof}
As a result of the previous propositions, we classified all pencils of cubics with only one base point. 

\section{Singular elements and non-isotriviality}\label{sec: singular}
Let $V\subset H^0(\pr^2,\cO_{\pr^2}(3))$ be any pencil of cubics satisfying all the conditions of \Cref{quest: main}. We now study in detail the singular elements of $|V|$. From this study we can conclude that $|V|$ is not isotrivial, i.e. that its smooth members are not mutually isomorphic.
\begin{proposition}\label{prop: nocusps}
The pencil $|V|$ has no cuspidal elements.
\end{proposition}
\begin{proof}
Let us first prove that, if there is a cuspidal cubic in $|V|$, the cusp cannot be $p$. 
Assume by contradiction that this is the case and call $D$ this cuspidal cubic. We have that $(D\cdot C)=9p$ and that $(D\cdot L)=3p$. 
In principle, we would like to use this information to draw a contradiction using \Cref{lem: namba}. However, $D$ is singular in $p$, so the Lemma does not apply straight away: we need to resolve the singularity of $D$ (see also \Cref{rem: singular-namba}). 
Let us call $\nu\colon X\to \pr^2$ the blow up of $\pr^2$ in $p$. Let $\widetilde{C}, \widetilde{D}, \widetilde{L}$ be the strict transforms of $C, D$ and $L$, respectively. Let $q$ be the point in the exceptional divisor $E$ corresponding to the slope of $L$, i.e. such that $q=\widetilde{L}\cdot E$. 

Let us verify that $(\widetilde{D}\cdot \widetilde{C})_q = 7$. Observe that  $(\widetilde{D}\cdot E)_q=2$, $(\widetilde{C}\cdot E)_q=1$. We have that 
\begin{align*}
9 &=(\nu^*(C)\cdot \nu^*(D))=((\widetilde{C}+E)\cdot (\widetilde{D}+2E))=\\
&=(\widetilde{D}\cdot \widetilde{C})+(E\cdot \widetilde{D})+2(E\cdot \widetilde{C})+2E^2=\\
&=(\widetilde{D}\cdot \widetilde{C})+2+2-2=(\widetilde{D}\cdot \widetilde{C})+2.
\end{align*}
We thus get a contradiction: indeed by \Cref{lem: namba}, as $\widetilde{D}$ is smooth, we should have 
\[
2=\min\{(\widetilde{D}\cdot \widetilde{C})_q, (\widetilde{D}\cdot E)_q\}\leq (\widetilde{C}\cdot E)_q=1.
\]

Let us now prove that there is no cuspidal elements in $V$ whose cusp $c$ is not $p$. 
Assume by contradiction that such an element $D$ exists. Observe that a cuspidal cubic has one smooth flex, which we call $f\in D$.
By similar arguments as above, we see that $f\not = p$. Indeed, if this is the case, let us resolve the singularity of $D$ by blowing up the cusp $c\in D$; call $\widetilde{C}, \widetilde{D}, \widetilde{L}$ be the strict transforms of $C, D$ and $L$ respectively, and let $\tilde p$ be the inverse image of $p$ in the blow up.
Clearly $\widetilde{D}$ is smooth and $(\widetilde{D}\cdot \widetilde{C})_{\tilde p}=(\widetilde{D}\cdot \widetilde{C})=9$, $(\widetilde{D}\cdot \widetilde{L})=(\widetilde{D}\cdot \widetilde{L})_{\tilde p} =3$, so by \Cref{lem: namba} again, 
\[
3=\min\{(\widetilde{D}\cdot \widetilde{C})_{\tilde p}, (\widetilde{D}\cdot  \widetilde{L})_{\tilde p}\}\leq (\widetilde{C}\cdot  \widetilde{L})_{\tilde p}=2,
\]
which gives a contradiction.

If we remove the cusp $c\in D$, we obtain a manifold homeomorphic to $\mathbb C$. 
Moreover, if we consider on $D\setminus\{c\}$ the sum defined as $d+d'=d''$ if and only if  $\cO_{D}(d'' )\cong \cO_D(d+d'-f)$, we have that  the group law on $\pic^0(D)$ induces an additive group law on $D\setminus\{c\}$ that makes it isomorphic to $\mathbb C$ with $f$ corresponding to $0$ (see for instance \cite[Appendix]{diller}).

Now, in $D$ we have that $\cO_D(9p)\cong \cO_D(9f)$, because $(D\cdot C)=9p$, $(D\cdot 3L)=9f$ and both $C$ and $3L$ belong to $|H^0(\pr^2,\cO_{\pr^2}(3))|$. So, we would have 
 in $D\setminus \{c\}=\mathbb C$ that $9p=0$, which would imply that $p=f$, a contradiction.
\end{proof}
\begin{remark}\label{rem: kol-cusp}
The second part of the proof, that excludes the existence of a cuspidal cubic with the cusp different from the base point $p$, can be proved also with Lemma 42 of \cite{Kollar}.
\end{remark}
\begin{proposition}\label{prop: singular fibres}
The pencil $|V|$ has 4 irreducible nodal cubics as singular elements, one of which has its node in the base locus $p$.
\end{proposition}
\begin{proof}
Consider the (topological) fibration obtained by removing the base point $p$ from $\pr^2$:
\[ \pi\colon \pr^2\setminus \{p\}\longrightarrow \pr^1.\]
The map $\pi$ has general fibres $F$ which are tori without a point, and $s$ fibres $F'$ which are nodal cubics minus one point. 
Note that in order to compute the topological Euler characteristic it is not relevant whether or not the removed point is the node because in any case we obtain a space homotopically equivalent to $S^1$.
The topological Euler characteristic of these fibres is $\chi_{top}(F)=-1$ and $\chi_{top}(F')=0$, respectively. 
Hence, the formula for the topological Euler characteristic for fibrations (see for instance \cite[Lemma VI.4]{beauville-libro}) tells us the following:
\begin{align*}
2=\chi_{top}(\pr^2\setminus \{p\})&=\chi_{top}(\pr^1)\chi_{top}(F)+s(\chi_{top}(F')-\chi_{top}(F))= -2+s.
\end{align*}
So, we obtain that $s=4$. 

Let us now prove that one of these fibres has node in $p$. Consider again the tangent line $L$ to $C$ at $p$. Any smooth member of the pencil has $L$ as a simple tangent at $p$. Indeed, if it exists $D$ smooth in $|V|$ such that $(D\cdot L)=3p$, then by \Cref{lem: namba} we have 
\[
3=\min\{(D\cdot C)_p, (D\cdot L)_p\}\leq (C\cdot L)_p=2,
\]
which gives a contradiction. 
Let us now consider the morphism $\varphi\colon |V|\to L$ that associates to any $D\in |V|$ the third point of intersection $(D\cdot L)= 2p+\varphi (D)$.
The map $\varphi$ is a morphism of projective varieties and it is injective, indeed if we had $\varphi(D)=\varphi(D')$ for $D\not =D'$ in $|V|$, then $\varphi(D)$ would be a base point different from $p$. So, $\varphi$ is an isomorphism (it is indeed a section of the pencil, as we verify in \Cref{prop: confronto}), and thus there exists an element $D'\in |V|$ such that $\varphi(D')=p$, i.e.  $(D'\cdot L)= 3p$. This curve $D'$ is necessarily singular by what observed above, and so it is a nodal cubic, with $p$ as the node and one of the branches tangent to $L$.
\end{proof}
\begin{remark}\label{rem: singular-namba}
From the argument of the above proposition, we see that there is an element $D'$ of the pencil $|V|$, singular in $p$, such that, by calling $L$ the tangent line to $C$ at $p$,
we have $(L\cdot D')=3p$,
and so 
\[\min\{(L,\cdot D')_p, (C,D')_p\}=3> 2=(C\cdot L)_p,\]
hence the formula of \Cref{lem: namba} does not hold.
This shows that the assumption that $D$ is smooth in $p$ in \Cref{lem: namba} is necessary.
\end{remark}
\begin{corollary}\label{cor: non-isotriviality}
The pencil $|V|$ is non-isotrivial.
\end{corollary}
\begin{proof}
As the pencil has (three) irreducible nodal members whose node is not the base point, they are fibres of the associated fibration, which is therefore non-isotrivial. 
Another way to prove the statement is to observe that the $j$-invariant of the nodal members of $|V|$ is $\infty$. 
Hence the modular map $J\colon \bP^1 \to \cM_{1,1} \cong \bP^1$ has poles at these fibres, so it necessarily is surjective.
\end{proof}

\begin{remark}\label{rem: isotriviality} 
It is important to stress that it can happen that some of the pencils described in \Cref{rem: esempio ovvio} can be isotrivial. Consider indeed the Fermat cubic $F(x,y,z)=y^2z-x^3+z^3$, and the pencil generated by $F$ and by $z^3$. Note that the line $z=0$ is tangent to $F=0$ at the point $(0,:1:0)$, so it is a pencil of the type describes in  \Cref{rem: esempio ovvio}. It is easy to see that this pencil is isotrivial; indeed it equation is of the form 
$y^2z-x^3+z^3(1+t)$, for $t\in \bC$, and so for any $t\not = -1 $ we have that if we choose a triple root $\zeta(t)$ of $t+1$, and a square root of $\zeta(t)$, say $\eta(t)$ we have that the projective transformation 
$(x:y:z)\mapsto(x:\eta(t)y:z/\zeta(t))$ sands the curve corresponding to $t$ to the original Fermat.
Otherwise one could also compute the $j$-invariant and see that it is the same for all smooth fibres.
The general member of any isotrivial pencil must have non-trivial automorphisms, as is the case of the Fermat, by the results of Serrano \cite{serrano}. See also \Cref{prop: flex}.
\end{remark}

\section{Comparison with a result of Gattazzo and other examples}\label{ssec: gattazzo}

As observed in \Cref{rem: sex}, \Cref{quest: main} can be rephrased for cubics as follows.  Given a smooth cubic $C$ and $p\in C$ we are looking for the existence of a different irreducible cubic such that the intersection of the two curves is $9p$.
This question was addressed in  \cite{Gat}.
We now recall some of the results contained in that paper, which is  written in Italian. Gattazzo's interesting construction also very naturally produces explicit equations for the pencils we are interested in.
\begin{definition}\label{def: type9}
Given a smooth cubic $C$, fix the origin of the group law in a flex $O$. The $72$ points of $9$-torsion which are not $3$-torsion are called {\em points of type $9$}.
\end{definition}
\begin{definition}\label{def: tang tri}
Let $C\subset \pr^2$ be a smooth cubic. A \emph{tangential triangle} for $C$ is a triple of points $p_0,p_1,p_2$ of $C$ such that 
$$(t_{p_0}C\cdot C)=2p_0+p_1, \qquad (t_{p_1}C\cdot C)=2p_1+p_2, \qquad (t_{p_2}C\cdot C)=2p_2+p_0,$$
where we denote by $t_{p}C$ the tangent line to $C$ at $p$.
\end{definition}

\begin{figure}[ht]
\includegraphics[width=140pt]{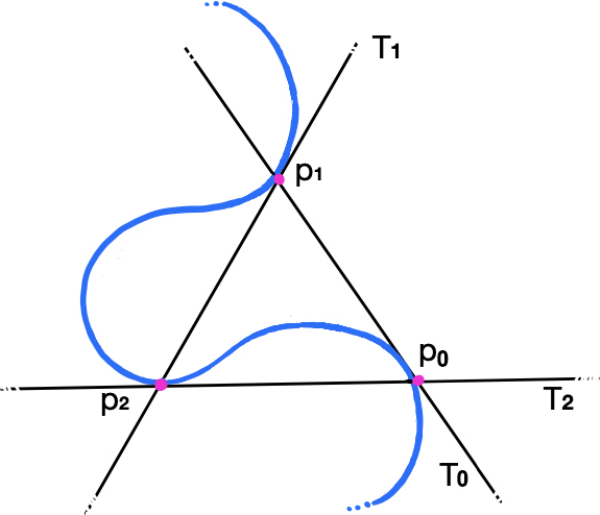}
  \caption{Gattazzo's tangential triangle}
  \label{fig: gattazzo}
\end{figure}

\begin{proposition}\label{prop: gatt}[Prop.~(2.1) of \cite{Gat}]
The point $q\in C$ is the vertex of a tangential triangle if and only if it is a point of type $9$.
\end{proposition}
\begin{proof}
Let $p_0, p_1,p_2$ be vertex of a tangential triangle, with associated lines $T_0,T_1,T_2$ as in \Cref{fig: gattazzo}. 
These point necessarily are not $3$ or $6$-torsion, by the definition of tangential triangle.  
Let $t_i\in \bC[x,y,z]/I(C)=A[C]$  be the linear equations of $T_i$, with $i=0,1,2$, where $I(C)$ is the ideal of the curve $C$, and $A[C]$ the coordinate ring of the affine cone associated with $C$. Now consider the 
meromorphic section ${(t_0^4t_2)}/{t_1^2}$ of $\cO_{\pr^2}(3)$.
By restricting the divisor of this section to  $C$ we have:
\begin{equation}\label{eq: gatt}
\divi\left(\frac{t_0^4t_2}{t_1^2}\right)=8p_0+4p_1+2p_2+p_0-(4p_1+2p_2)=9p_0.
\end{equation}
So, $h:=(t_0^4t_2)/t_1^2$ is a meromorphic section whose divisor $\divi(h)$ is effective. 
This implies that  $h\in A[C]_{\mathscr{P}}$ for any  $\mathscr{P}$ prime ideal in $A[C]$ of height $1$. By the structure theorem of integrally closed domains we have that 
\[A[C]=\bigcap_{\mathscr{P},\,\, ht(\mathscr{P})=1}A[C]_{\mathscr{P}}.\]
So, we have that $h\in A[C]$. 
Consider a polynomial $\mathscr{H}\in \bC[x,y,z]$ whose image in $A[C]$ via the canonical projection is $h$. 
It defines a  curve $D_3$ in $\pr^2$ of degree $\deg\divi(h)/3=3$ such that $(C\cdot D_3)=\divi(h)=9p_0$.
It is easy to prove that $D_3$ is irreducible. 
Conversely, by starting from a point $p$ of order $9$, and taking the tangent at the third point of intersection with $C$, and performing the same construction twice, one gets back to $p$ obtaining a tangential triangle.
\end{proof}
In \cite{Gat} is then proved that there are precisely 24 tangential triangles, with vertices all distinct, and this gives precisely the 72 points of type $9$.

\begin{example}\label{ex: Gatt}
Consider the cubic $C$ given by the polynomial
\[
xy^2+yz^2+zx^2+3bxyz, \mbox{ with }b^3\not=-1.
\]
It is easy to prove that the three coordinate points $p_1=(1:0:0)$, $p_2=(0:1:0)$, $p_3=(0:0:1)$ are the vertices of a tangential triangle for $C$.
Then we can directly write down the equation of the cubic $C_3$ in \Cref{prop: gatt} starting from the element in the fraction field whose divisor is $9p$ as in \eqref{eq: gatt}. 
Recall that $t_1=z$, $t_2=x, t_3=y$, and so we obtain 
\begin{equation*}
\begin{split}
\frac{y^4x}{z^2}&=\frac{y^2(-yz^2-zx^2-3bxyz)}{z^2}=\frac{y^2(-yz-x^2-3bxy)}{z}=-y^3-\frac{xy^2(x+3by)}{z}=\\
&=-y^3+\frac{(yz^2+zx^2+3bxyz)(x+3by)}{z}=-y^3+(yz+x^2+3bxy)(x+3by),
\end{split}
\end{equation*}
where the equalities are mod $I(C)$. 
So, in particular by choosing $b=0,1$ we get
\[ xy^2+yz^2+zx^2,  \quad x^3+xyz-y^3\]
and 
\[ xy^2+yz^2+zx^2+3xyz,  \quad -y^3+(yz+x^2+3xy)(x+3y).\]
These give two examples of pencils of cubics satisfying all the conditions of \Cref{quest: main}.
It is easy to check that the general members are smooth cubics.
\end{example}


\section{Uniqueness up to projective transformations}\label{ssec: uniqueness}
\begin{proposition}\label{prop: uniqueness}
Let us consider any pencil $|V|$ satisfying $(C1)$, $(C2)$ and $(C3)$. Then, there is a projective transformation of $\bP^2$ sending $|V|$ to the pencil of \Cref{ex: Gatt} generated by 
\begin{equation}\label{eq: pencil}
xy^2+yz^2+zx^2\quad  \mbox{ and }\quad x^3+xyz-y^3.
\end{equation}
\end{proposition}
\begin{proof}
Let $p_0$ be the base point of $|V|$. As $|V|$ is non-isotrivial by \Cref{cor: non-isotriviality}, there is a smooth member $C_0$ of $|V|$ that has the same $j$-invariant as the Klein cubic (and so the same as the Fermat cubic). Now, call $p_1$ and $p_2$ the other vertices of the tangential triangle on $C_0$ associated to $p_0$ (see \Cref{prop: gatt}).
Consider the projectivity of $\bP^2$ sending $p_0\mapsto (1:0:0)$, $p_1\mapsto (0:1:0)$ and $p_2\mapsto (0:0:1)$. 
Then, by imposing the tangency conditions, we have that an equation for the image of $C_0$ is necessarily of the form 
\[axy^2+byz^2+czx^2+dxyz=0,\]
with $a,b,c\not=0$. 
Now, if we apply the projectivity 
\[
\left\{
\begin{matrix}
x\mapsto \alpha x\\
y\mapsto \beta y\\
z\mapsto \gamma z\\
\end{matrix}
\right.
\]
where $\gamma$ is a fifth root of $a^2/(cb^4)$, $\alpha=(b^2\gamma^4)/a$ and $\beta =1/(b\gamma^2)$, we obtain an equation of the form
\[xy^2+yz^2+zx^2+\zeta xyz=0,\]
for some $\zeta \in \bC$. Now recalling that by our initial assumption $C_0$ has the same $j$-invariant as the Klein cubic, we have necessarily $\zeta=0$.

Now we just observe that, by \Cref{prop: construction} and \Cref{prop: unicita'}, there is only one pencil containing the Klein cubic $xy^2+yz^2+zx^2=0$, with one base point in $(1:0:0)$ satisfying (C1), (C2) and (C3), 
which is exactly \eqref{eq: pencil}.
\end{proof}

Throughout the rest of the paper, we denote by $|\cV|$ the pencil with generators in \eqref{eq: pencil}. Every pencil satisfying (C1), (C2) and (C3) admits a projective transformation sending it to $|\cV|$.


\begin{remark} \label{rem:kollardue}
An alternative proof of the uniqueness up to projective transformations follows from \cite[Example 46]{Kollar}: starting from a nodal cubic $C_0$ at $p$, there exists a unique pencil $|V|$ satisfying $(C1)$, $(C2)$ and $(C3)$ with $p$ as base point with all the members tangent to $C_0$ on a fixed branch at $p$. 

As suggested by J\'anos Koll\'ar in a private communication, this can be used to give a classification of this kind of pencils up to projective transformation over any field $k$, possibly of characteristic different from $2, 3$. 
Indeed, up to projective transformations, a nodal cubic on $\bP^2_k$ with tangent lines at the node defined over $k$ has equation 
\begin{equation} \label{eq:formacanonica}
x^3+\lambda y^3-xyz=0,
\end{equation}
where $\lambda$ is in $k^\ast$. Note that the line passing through all the nonsingular flexes over an extension of $k$ is defined on the field $k$ and we obtain Equation \eqref{eq:formacanonica} by imposing that such line has equation $z=0$.
Two of these equations associated to $\lambda$ and $\lambda'$ gives the same cubic if and only if $\lambda/\lambda'$ belongs to $(k^\ast)^3$.
As a consequence the pencils are classified by $k^\ast / (k^\ast)^3$.
\end{remark}

\section{Comparison with Beauville's pencils}\label{ssec: beauville}

From \Cref{prop: singular fibres} and \Cref{cor: non-isotriviality}, we deduce that, after blowing up 9 points infinitely close to $p$, the pencil $|\cV|$ gives rise to a semistable family of curves of genus $1$ over $\pr^1$ with exactly $4$ singular fibres. 
This is the minimal possible number of singular fibres of non-isotrivial semistable fibrations, as proved in \cite{beau1}. 

The fibration induced by $|\cV|$ has thus to be isomorphic to one of the $6$ families classified by Beauville in \cite{beauville}, that we list in \Cref{fig: beauville}.
  \bigskip
\begin{table}[ht]
    \begin{tabular}{l|c|c}\label{tab: beauville}
       & Equation & Number of components\\ && of singular fibres \\
      \hline
      $(F1)$ & $x^3+y^3+z^3+txyz=0$& $3,3,3,3$\\
      $(F2)$ & $x(x^2+z^2+2zy)+t(x^2-y^2)=0$ & $4,4,2,2$\\
        $(F3)$ & $x(x-z)(y-z)+tzy(x-y)=0$& $5,5,1,1$\\
                $(F4)$ &$(x+y)(y+z)(z+x)+txyz=0$ &$6,3,2,1$ \\
          $(F5)$ & $(x+y)(xy-z^2)+txyz=0$& $8,2,1,1$\\
            $(F6)$ & $x^2y+y^2z+z^2x+txyz=0$& $9,1,1,1$\\
    \end{tabular}
    \caption{Beauville's table of fibrations if genus $1$ with $4$ singular fibres over $\pr^1$.}
    \label{fig: beauville}
\end{table}

The resolution of the base point of $|\cV|$ is obtained by blowing up 9 times $p$ and points infinitely close to $p$, so it does not affect the $3$ singular fibres whose node is not $p$, nor the smooth fibres. 
By looking at the number of components of the fibres, the only possible family is (F6).
Note moreover that (F6), with $t=0$, is a special case of the first equation in Gattazzo's pencil in \Cref{ex: Gatt}. But the pencil of Beauville has $3$ base points: indeed the second element $xyz=0$ of Beauville's pencil is precisely the equation of a tangential triangle to the first curve.

We now explicitly verify that one of the singular fibres of the fibration induced by $|\cV|$ is a ring composed of $9$ rational components; we then exhibit an explicit Cremona transformation of order $4$ sending $|\cV|$ into the family (F6).

\begin{proposition}\label{prop: confronto}
Let $\nu \colon X\to \pr^2$ be the sequence of blow ups that provide a resolution of the base point $p$ of $|\cV|$, and let $f\colon X\to \pr^1$ be the induced fibration. The singular fibres of $f$ are three irreducible nodal curves of arithmetic genus $1$ and one curve of arithmetic genus $1$ made of a ring of $9$ rational $(-2)$-curves. The family $f$ is isomorphic to family (F6) of \Cref{fig: beauville}.
\end{proposition}
\begin{proof}
\begin{figure}[ht]
\includegraphics[width=\textwidth]{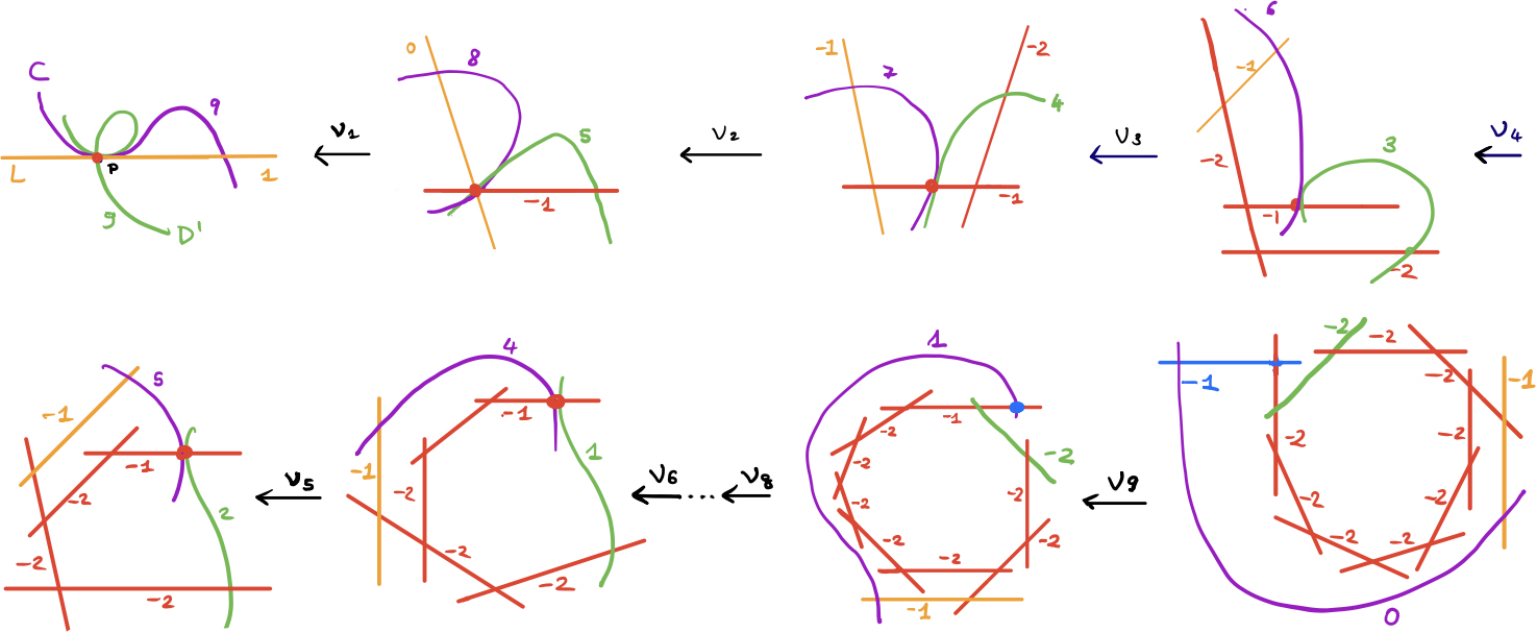}
  \caption{ Blowups resolving the base points of $|V|$. {\color[RGB]{136,48,187}$C$ is the smooth cubic in the  pencil}, {\color[RGB]{255,180,0}$L$ is the tangent line to $C$ at $p$} and {\color[RGB]{100,160,70}$D'$ is the nodal cubic with node in $p$}.}
  \label{fig: scoppiamenti-1}
\end{figure}
Proof by drawing it shall be.
In \Cref{fig: scoppiamenti-1}  we follow the nine blow ups, called $\nu_1,\ldots \nu_9$. 
We call $X$ the final surface and $\nu\colon X\to \pr^2$ the composition $\nu_9\circ \nu_8\circ \ldots \circ \nu_1$.
The curves represented in \Cref{fig: scoppiamenti-1} are the original curve $C$, the nodal curve with node in $p$, which we call as above $D'$, and the tangent line $L$ to $C$ at $p$. 
We have labelled each curve with its self-intersection. 
At the end of the process the tangent line $L$ becomes a section of the fibration, as it should be by the remarks made in \Cref{prop: singular fibres}. 
Also the last exceptional curve, coloured in blue in the picture, is a section of the fibration $f$. 
The singular fibre obtained at the end is  a closed chain made of the $8$ first exceptional curves (in red in the picture) plus the strict transform of $D'$. All these components are $(-2)$-curves. 

There are a couple of things worth verifying. First of all, we see in the picture that $\nu_8$ separates the inverse image of $C$ and $D'$. This comes from the fact that after $\nu_1$ the intersection of the strict transform of $D'$ and $C$ drops from $9$ to $7$. Indeed, let us call $\widetilde {D'_1}$ and $\widetilde {C}_1$ these two curves. We have 
\[\nu_1^*(C)=\widetilde {C_1}+E_1,\quad \nu_1^*(D')=\widetilde {D'_1}+2E_1,\]
where $E_1$ is the first exceptional divisor.
So we obtain 
\begin{align*}
&9=(\nu_1^*(C)\cdot \nu_1^*(D'))=(\widetilde {C}_1+E_1\cdot \widetilde {D'}_1+2E_1)=\\
&=(\widetilde {C}_1\cdot \widetilde {D'}_1)+2(\widetilde {C}_1\cdot E_1)+(\widetilde {D'}_1\cdot E_1)+ 2E_1^2=(\widetilde {C}_1\cdot \widetilde {D'}_1)+2.\end{align*}
Let us also verify that the self-intersection of $\widetilde{D'}$ is 5. Indeed we have 
\[ 9=(\nu_1^*(D'))^2=(\widetilde {D'_1})^2+4\widetilde {D'_1}E_1+4E_1^2=(\widetilde {D'_1})^2+8-4,\]
which gives $(\widetilde {D'_1})^2=5$ as stated.

Eventually, the last statement of tie proposition derives from Beauville's classification Theorem in \cite{beauville}, where he proves that the $6$ families in \Cref{fig: beauville} are unique up to isomorphism of fibrations.
\end{proof}
\begin{remark}\label{rem: kollar}
In the context of Enriques surfaces, the correspondence between \eqref{eq: pencil} and (F6) is described in \cite[Remark 8.9.31]{dk2}, see also \cite[(3.7)]{kondo}.
Similar pictures as \Cref{fig: scoppiamenti-1}  have been obtained in \cite{yos}, \cite[Sec.6]{ore} and in \cite[Sec.2.1]{MS}.
\end{remark}
We now exhibit an explicit Cremona transformation of $\pr^2$ sending (F6) to $|\cV|$, up to projective transformations. We firstly blow up the indeterminacy locus of (F6). 
The equation (F6) is given by the combination of $xyz=0$, three lines in general position, and $x^2y+y^2z+z^2x=0$, the Klein cubic that meets the lines in the three points of intersections $(1:0:0)$, $(0:1:0)$, and $(0:0:1)$, and these are the vertices of a tangential triangle as already observed in the beginning of this section.
\begin{figure}[ht]
\includegraphics[width=\textwidth]{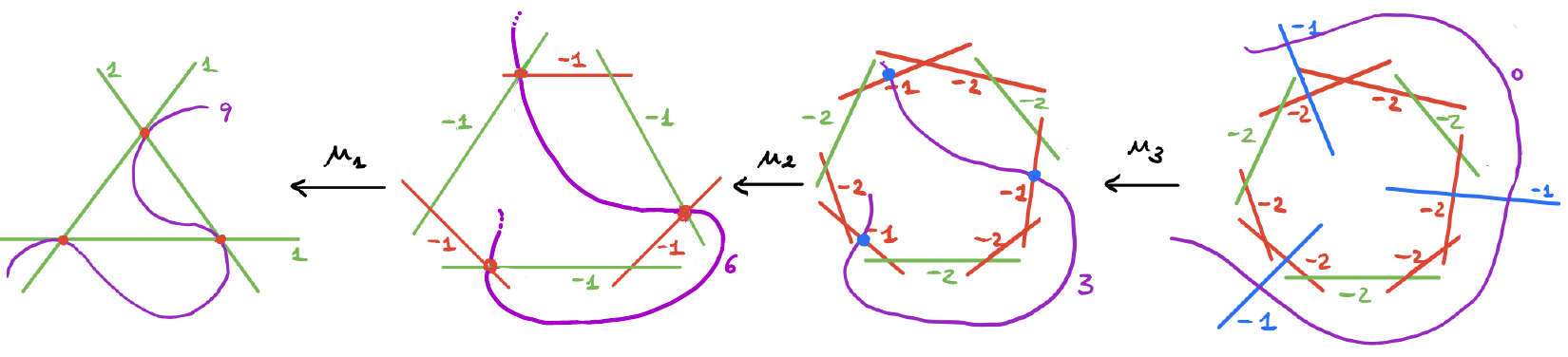}
  \caption{Blow ups resolving the base points of (F6) in \Cref{fig: beauville}.}
  \label{fig: scoppiamenti-2}
\end{figure}

In \Cref{fig: scoppiamenti-2}, $\mu_1$ consists of three (simultaneous) blow ups in the three coordinate points, $\mu_2$ of  three (simultaneous) blow ups in three points each in one of the exceptional divisors of $\mu_1$, and $\mu_3$ of three (simultaneous) blow ups in three points each in one of the exceptional divisors of $\mu_2$. Let us call $\mu\colon Y\to \pr^2$ this chain of blowups, and call $g\colon Y\to \pr^1$ the fibration obtained.
As we see from the figure, we obtain the same singular fibre as in \Cref{fig: scoppiamenti-1}, and three sections of $g$ with the last three blow ups (in blue in the picture).
By Beauville's result \cite{beauville}, there is an isomorphism of fibrations $\eta$ making the following diagram commute. 
\[
\xymatrix{
X\ar_f[dr]\ar^{\eta}[rr]&&Y\ar^g[dl]\\
&\pr^1&}
\]

\begin{figure}[ht]
\begin{center}
\includegraphics[width=125pt]{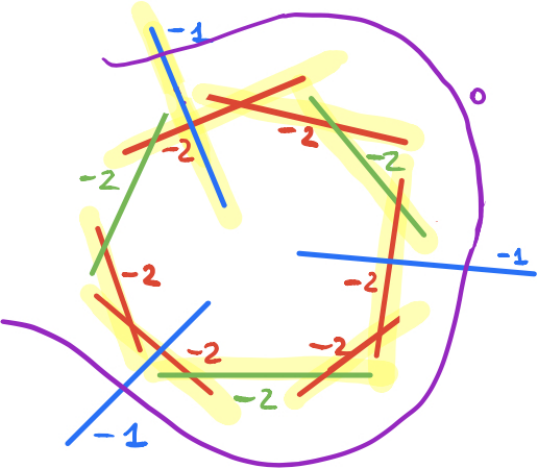}
\caption{The chain of rational curves to be contracted in $Y$ to get back $|\cV|$.}\label{fig: contrazioni-1}
\end{center}
\end{figure}

\begin{figure}[ht]
\begin{center}
\includegraphics[width=380pt]{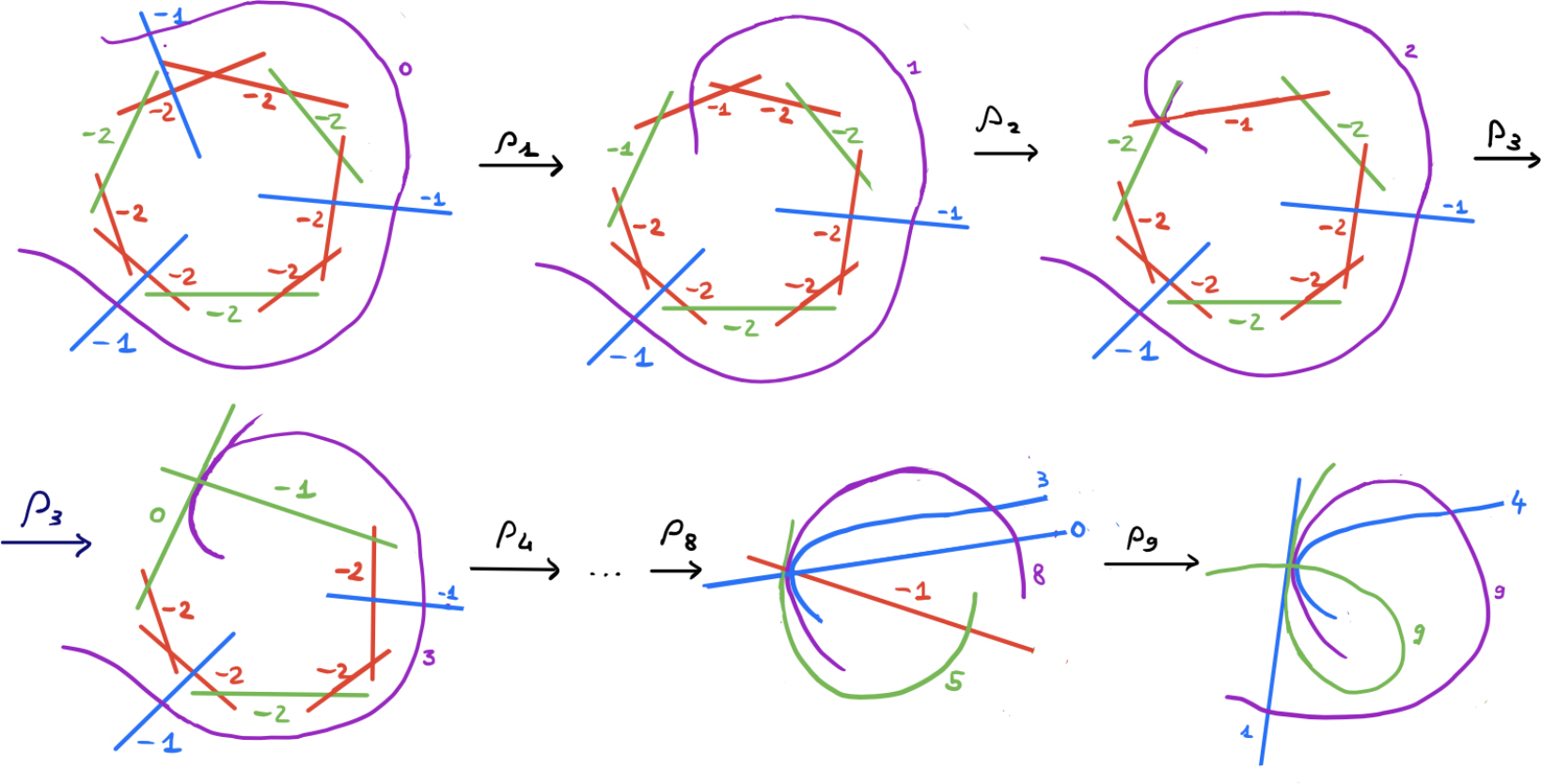}
\caption{Blow ups $\rho\colon Y\to \pr^2$ that send pencil (F6) in \Cref{fig: beauville} to $|\cV|$.}
\label{fig: contrazioni-3}
\end{center}
\end{figure}

We can be more explicit:
it is now clear how to perform an explicit blowing down of $Y$ that induces, together with $\mu$, a Cremona transformation $\phi\colon \pr^2 \dasharrow \pr^2$ sending Beauville's pencil to $|\cV|$. Choose the chain of rational curves in $Y$ highlighted in 
\Cref{fig: contrazioni-1}.
This is a chain of rational curves where the first one (the chosen section) is a $(-1)$-curve, and all the others are $(-2)$-curves. So, by blowing down the chain in order, starting from the  $(-1)$-curve, we go back to a rational surface with Picard number $1$, so to $\pr^2$. If we follow Beauville's pencil through these blow downs, we see that we  obtain a pencil whose general elements are smooth and intersect in one point, with order of tangency $9$. The $(-2)$-curve that we did not blow down in the singular fibre of Beauville's family  has become a nodal cubic having a node in the base point, as in $|\cV|$ (see \Cref{fig: contrazioni-3}).

We call $\rho\colon Y\to \pr^2$ this chain of blow ups.
The composition $\phi:=\rho\circ \mu^{-1}$  is thus a Cremona transformation  sending (F6)  to $|\cV|$. 
\[
\xymatrix{
&Y\ar_{\mu}[dl]\ar[dr]^{\rho}&\\
\pr^2\ar^{\phi}@{-->}[rr]&&\pr^2}
\]
\begin{figure}[ht]
\begin{remark}
We can moreover observe the following: in the definition of $\mu\colon Y\to \pr^2$ we have blown up simultaneously the three points (and infinitely close points to them), but indeed, we could also change the order of blowing up, by blowing up first one point and 2 infinitely closed points and then the other and last the third. If we choose as final the point that gives rise to the chosen section, we have that the last three blow ups of $\mu$ inverse to to the last three blow down of $\rho$, so we can skip them, and contract only a chain of $6$ curves in the  blow up of $\pr^2$ in 6 points, as described by \Cref{fig: contrazioni-2}. 
\end{remark}
\begin{center}
\includegraphics[width=130pt]{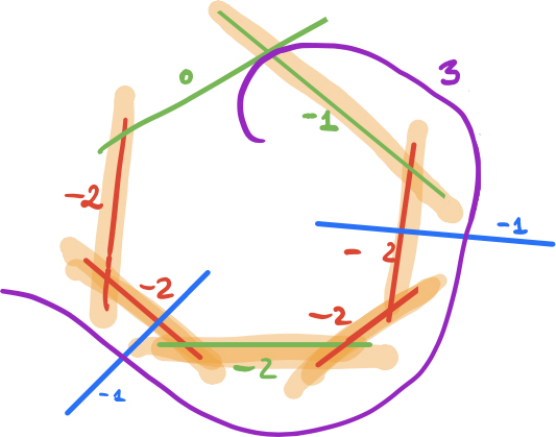}
\end{center}
\caption{Chain of contractions of only $6$ curves to get back $|\cV|$.}
\label{fig: contrazioni-2}
\end{figure}

So we obtain the following
\begin{proposition}\label{prop: cremona 4}
There is a Cremona transformation $\phi\colon \pr^2\dasharrow \pr^2$ of degree $4$  sending $|\cV|$ to  Beauville's pencil (F6).
\end{proposition}
\begin{proof}
The only thing to prove is the statement on the degree of $\phi$. This is equivalent to the fact that a general line in $\bP^2$ is sent by $\phi$ in a curve of degree $4$. So let $\ell$ be a general line. Its strict transform $\widetilde \ell$ via $\rho$ is still a curve in $Y$ with self-intersection $1$. Now, $\widetilde \ell$ intersects transversally $E_3$ in one point  and $E_6$ in two points, so via the blow downs of the right hand side of \Cref{fig: contrazioni-1} we obtain that the image curve has self-intersection $16$. Hence the Cremona $\phi$ has degree $4$.
\end{proof}

It is now clear how to produce the inverse Cremona sending $|\cV|$ to the family (F6).
In order to do this explicitly, we need to find the three sections that has to be blown down first. We already have two sections of the fibration $f\colon X\to \pr^1$: 
one is given, coherently with the observation in the proof of \Cref{prop: singular fibres}, by the strict transform of the tangent $L$ to $C$ at $p$; the second section is the exceptional curve of  $\nu_9$. 

We see in \Cref{fig: contrazioni-3} what happens to the other section of Beauville's family $g\colon Y\to \pr^1$: it becomes a rational curve that has tangency $5$ with $C$ at $p$, and self-intersection $4$. So, it is the unique osculating conic $Q$ to $C$ at $p$ (see \Cref{rem: sex}).  It is easy, via the very same reasoning made in \Cref{sec: singular}, to see that $Q$ need to be smooth.
Its strict transform in $X$ is the third section of $f$. 

\begin{figure}[ht]
\begin{center}
\includegraphics[width=140pt]{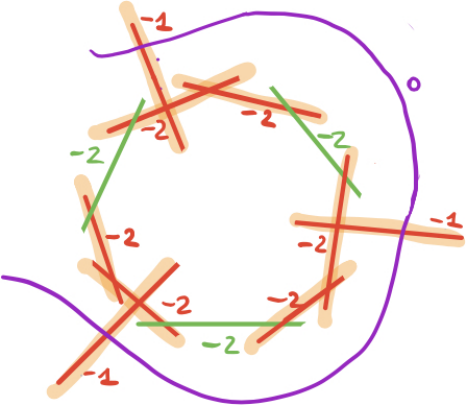}
\end{center}
\caption{The three chains of rational curves to be contracted to obtain the pencil of Beauville.}
\label{fig: contrazioni-4}
\end{figure}

In order to obtain (F6) we thus simply blow down the three chains of rational curves in \Cref{fig: contrazioni-4}, starting from the three sections.

\section{Pencils with one base point which is a flex}
With the same arguments as in \Cref{prop: nocusps,prop: singular fibres} we can prove the following:
\begin{proposition}\label{prop: flex}
Let $|W|$ be a pencil of plane cubics satisfying $(C1)$ and $(C2)$ whose base point $p$ is a flex for its smooth members. Let $L$ be the tangent line at $p$ to any smooth member of $|W|$. Then the following hold:
\begin{enumerate}
\item[(i)] the pencil $|W|$ has no other reducible member apart from $3L$;
\item[(ii)] the nodes and cusps of the singular elements are not $p$;
\item[(iii)] the possibilities for its singular members are:
\begin{enumerate}
\item $3L$ and two nodal cubics if  and only if $|W|$ is non-isotrivial;
\item $3L$ and one cuspidal cubic if and only if $|W|$ is isotrivial.
\end{enumerate}
\end{enumerate}
\end{proposition}
\begin{proof}
The proof for (i) is identical to the one in \Cref{prop: construction}, so we skip it. 

As for (ii), if $p$ were a singular point of another member of the pencil then all elements would be singular at $p$.

To prove point (iii), let us consider the topological fibration obtained by removing $p$ from $\pr^2$:
\[ \pi\colon \pr^2\setminus \{p\}\longrightarrow \pr^1.\]
Call $F$  the general fibres of $\pi$, call $F'$ the nodal cubics without the point $p$, and $F''$ the cuspidal cubics without the point $p$. 
As in \Cref{prop: singular fibres}, we have $\chi_{top}(F)=-1$, $\chi_{top}(F')=0$. On the other hand, a plane cubic with a cusp is homeomorphic to $\pr^1_\bC$, and so $\chi_{top}(F'')=1$. 
Then, consider the formula for the topological Euler characteristic for fibrations
(\cite[Lemma VI.4]{beauville-libro}):
\begin{align*}
2&=\chi_{top}(\pr^2\setminus \{p\})=\chi_{top}(\pr^1)\chi_{top}(F)
+(\chi_{top}(L\setminus\{p\})-\chi_{top}(F))+\\
&\quad +s'(\chi_{top}(F')-\chi_{top}(F))+s''(\chi_{top}(F'')-\chi_{top}(F))=\\
&=-2+2+s'(\chi_{top}(F')-\chi_{top}(F))+s''(\chi_{top}(F'')-\chi_{top}(F)),
\end{align*}
where $s'$ (resp. $s''$) is the number of fibres of type $F'$ (resp. $F''$).
So we see that the only possibilities are $(a)$ and $(b)$. 

The first case is non-isotrivial as in \Cref{cor: non-isotriviality}. In the second case the pencil is isotrivial as both singular fibres have finite $j$-invariant, so the whole modular map is constant.
For the second case, one could also derive isotriviality by observing that  the associated fibration has two singular fibres, while the minimum number of singular fibres of any non-isotrivial genus $g\geq 1$ fibration is $3$ by \cite[Proposition 1]{beau1}.
\end{proof}

\begin{remark}
There exist pencils as in (a) and (b) of the previous proposition. Both are constructed as in \Cref{rem: esempio ovvio}. Case (a) is given by starting with a nodal cubic, and case (b) is as in \Cref{rem: isotriviality}.
\end{remark}

\section{Comparison with Miranda-Persson classification} \label{sec: Miranda-Persson}

We have constructed pencils of cubics satisfying (C1) and (C2).
As above, we call $|\cV|$ the pencil satisfying also (C3), and we will call $|W|$ (respectively $|W'|$) two pencils as in $(a)$ (respectively $(b)$) of \Cref{prop: flex}. Recall that  $|W|$ is non-isotrivial and $|W'|$ is isotrivial. 
All the fibrations associated have genus $1$, and we will see that they are rational elliptic fibrations (i.e. they admit at least one section). 
The possible fibres of elliptic fibrations have been classified by Kodaira \cite{K}, and Persson classified the fibrations themselves in \cite{per}.
Let $f\colon X\to B$ an elliptic fibration with a distinguished section (sometimes this is referred to as Jacobian fibration). 
Recall that the set $\Phi(X)$ of sections of $f$ have the structure of an abelian group with the distinguished section as the zero element, usually called the Mordell-Weil group.

\begin{definition}[\cite{MP}]\label{def: reef}
An elliptic fibration $f\colon X\to B$ is called {\em extremal} if 
\begin{itemize}
\item[(i)] $\rho(X)=h^{1,1}(X)$ (i.e. the Picard rank is maximal);
\item[(ii)] the group of sections $\Phi(X)$ has rank $0$ (i.e. it is torsion).
\end{itemize}
\end{definition}
Let us call $f, g, g'\colon X\to \bP^1$, the fibrations associated to $|\cV|$, $|W|$ and $|W'|$, respectively and $X$ is the blow up of $\pr^2$ in 9 infinitely closed points. 
\begin{proposition}
The fibrations $f$, $g$ and $g'$ are all rational extremal elliptic fibrations. They correspond in Miranda-Persson's list in \cite[Theorem 4.1]{MP} to the elements $X_{1119}$, $X_{211}$ and $X_{22}$, respectively.
\end{proposition}
\begin{proof}
We need to understand the fibres of $g$ and $g'$. Recall that $L$ is the tangent line at $p$ to any smooth member of the pencils $|W|$ and $|W'|$. We need to understand how the element $3L$ transforms after the successive blow ups. 
\begin{figure}[ht]
\begin{center}
\includegraphics[width=350pt]{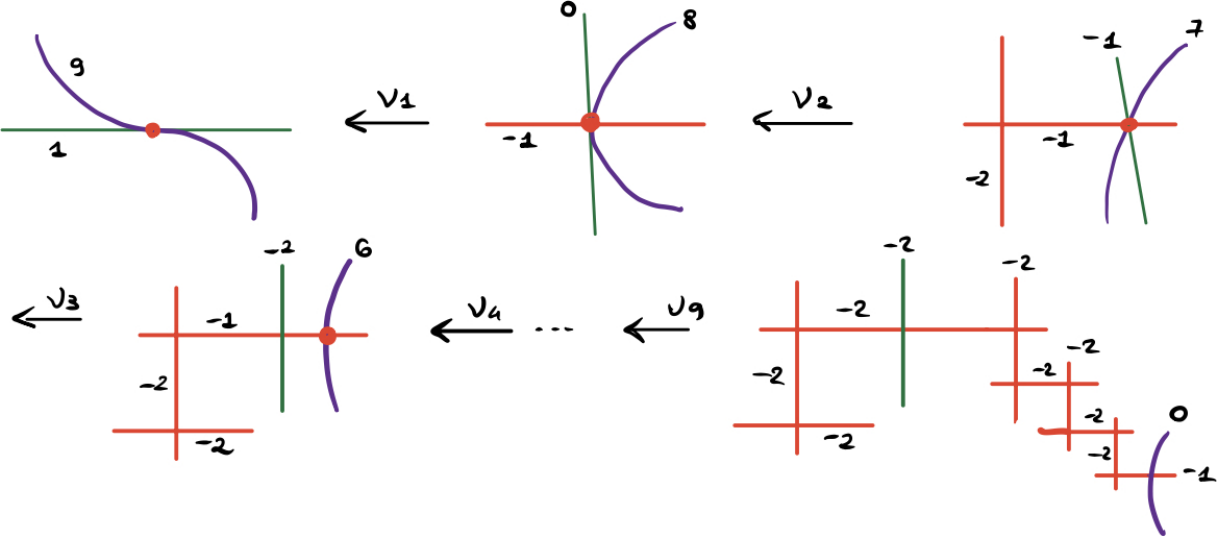}
\end{center}
\caption{The consecutive blow ups and the configuration of the exceptional divisors keeping track of the element $3L$.}
\label{fig: fibra-1}
\end{figure}
We see in \Cref{fig: fibra-1} the consecutive blow ups and the configuration of the exceptional divisors. We keep track of the self-intersection at any step. 
Call $\nu\colon X\to \bP^2$ the chain of blow ups, and call $E_i$ the exceptional divisors. 
The total transform of $3L$ is:
\[
\nu^*(3L)=3E_1+6E_2+9(E_3+\ldots E_9)+3\widetilde L,
\]
where $\widetilde L$ is the strict transform of $L$. So, the singular fibre corresponding to $3L$ is:
\[
\nu^*(3L)-\sum_{i=1}^9iE_i=3\widetilde L + 2E_1+4E_2+6E_3+5E_4+4E_5+3E_6+2E_7+E_8,
\]
corresponding to \Cref{fig: fibra-2}. This is a fibre of type $II^*$ in the classification of Kodaira.
\begin{figure}[ht]
\begin{center}
\includegraphics[width=160pt]{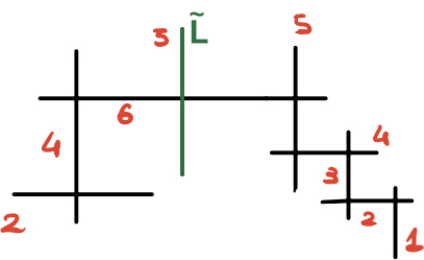}
\end{center}
\caption{The singular fibre corresponding to $3L$.}
\label{fig: fibra-2}
\end{figure}
We see that in this case the last exceptional divisor $E_9$ is again a section of  both $g$ and $g'$, so that they are elliptic rational fibrations. 

Now we prove that all our fibrations are extremal. Recall the formula of Shioda-Tate (see \cite[Corollary 2.4]{MP}), holding for any elliptic fibration $h\colon X\to B$:
\begin{equation}\label{eq: ST}
\rho(X)=2+\rk\Phi(X)+\sum_{F'}r(F'),
\end{equation}
where the sum goes over all the reducible fibres $F'$ of $h$, and $r(F')$ is the number of components of the support of $F'$ minus $1$.
In our case we have that $\rho(X)=10$, and in all three fibrations we have only one non-reduced fibre with $9$ components. As a consequence formula \eqref{eq: ST} becomes, for all the fibrations $f,g,g'$:
\[
10=\rho(X)=2+\rk\Phi(X)+8.
\]
Thus the group $\Phi(X)$ is torsion, and the fibration is extremal, as wanted.
In order to find the fibrations in the list of Miranda and Persson, we observe what are the types of the singular fibres using the Kodaira notation for elliptic fibrations:
\begin{itemize}
\item $f$, corresponding to $|\cV|$, has fibres $I_9, I^3_1$;
\item $g$, corresponding to $|W|$, has fibres $II^*, I^2_1$;
\item $g'$, corresponding to $|W'|$, has fibres $II^*, II$.
\end{itemize}
These configurations of fibres identify uniquely the elements $X_{1119}$, $X_{211}$ and $X_{22}$, respectively. 
\end{proof}
\begin{remark}
We can see the isotriviality of the fibration $g'$ also in the list of \cite[Theorem 4.1]{MP}, because the modular map $J\colon \bP^1\to \overline{\mathcal{M}}_{1,1}\cong \bP^1$ has degree $0$.
On the other hand, the modular map for the fibration $f$ corresponding to $|\cV|$ has degree $12$, as proved in  {\it loc.~cit.} 
In this case, the three sections that compose the Mordell-Weil group $\Phi(X)$ identify a tangential triangle.
Given a smooth plane cubic $C$ and fixing a point $p$ of type $9$ on $C$, there exists a fibre $F_t$ of $f$ with the $3$ points $\Phi(X)\cap F_t$ corresponding to the tangential triangle associated to $p$. With this choice $[(C,p)]\in \mathcal M_{1,1}$ is the image of the modular map corresponding to the fibre $F_t$. 
It seems that the preimage $J^{-1}([(C,p)]$ should identify a subset of 12 tangential triangles on $C$. 
So, there seems to be a natural partition of the tangential triangles into two subsets of cardinality 12 each, that might be further investigated.
\end{remark}

\bigskip

\noindent Riccardo Moschetti,\\Dipartimento di Matematica, Universit\`a di Pavia, Italy.\\
E-mail: \textsl {riccardo.moschetti@unipv.it}.
\medskip

\noindent Gian Pietro Pirola,\\Dipartimento di Matematica, Universit\`a di Pavia, Italy.\\
E-mail: \textsl {gianpietro.pirola@unipv.it}.
\medskip

\noindent Lidia Stoppino,\\Dipartimento di Matematica, Universit\`a di Pavia, Italy.\\
E-mail: \textsl {lidia.stoppino@unipv.it}.

\end{document}